\documentclass[12pt]{amsart}
\usepackage{amssymb}
\usepackage{amsmath}
\usepackage{amsthm}
\numberwithin{equation}{section}
\theoremstyle{plain}
\newtheorem{theorem}{Theorem}[section]
\newtheorem{proposition}[theorem]{Proposition}
\newtheorem{lemma}[theorem]{Lemma}
\newcommand{\C}{\ensuremath{\mathbb{C}}}

\def\cal{\mathcal}

\def\L{\Lambda}

\def\p{\partial}

\def\W{{\cal W}}
\begin{document}



\vspace*{-1cm}
\hbox{ }
{\hspace*{\fill} Mannheimer Manuskripte 260}

{\hspace*{\fill} math/0102169}

\vspace*{2cm}

\title{Almost K\"ahler deformation quantization}

\author[A.V. Karabegov]{Alexander V. Karabegov}
\address[Alexander V. Karabegov]{Theory Division, Yerevan 
Physics
Institute, Alikhanyan bros. 2, Yerevan 375036, Armenia}
\email{karabeg@uniphi.yerphi.am}
\author[M. Schlichenmaier]{Martin Schlichenmaier}
\address[Martin Schlichenmaier]{Department of Mathematics and 
  Computer Science, University of Mannheim, A5 \\
         D-68131 Mannheim \\
         Germany}
\email{schlichenmaier@math.uni-mannheim.de}
\begin{abstract}
We use a natural affine connection with nontrivial torsion on 
an
arbitrary almost-K\"ahler manifold which respects the 
almost-K\"ahler
structure to construct a Fedosov-type deformation quantization 
on this
manifold.
\end{abstract}
\subjclass{Primary: 53D55; Secondary:  53D15, 81S10, 53D50, 53C15}
\keywords{deformation quantization, star product, 
almost-K\"ahler manifold.
}

\date{February 20, 2001, rev. June 18, 2001}
\maketitle

\section{Introduction}\label{S:intro}
Let $(M,\{\cdot,\cdot\})$ be a Poisson manifold. Denote by 
$C^\infty(M)[[\nu]]$ the space of formal series in 
$\nu$ with coefficients from $C^\infty(M)$.
Deformation quantization on $(M,\{\cdot,\cdot\})$, as defined 
in \cite{BFFLS}, is an associative algebra structure on 
$C^\infty(M)[[\nu]]$ with the $\nu$-linear and $\nu$-adically 
continuous product $\ast$ (named a star-product) given on 
$f,g\in C^\infty(M)$ by the formula
\begin{equation}\label{E:star}
 f \ast g = \sum_{r=0}^\infty \nu^r  C_r(f,g),
\end{equation}
where $C_r, r\geq 0$, are bilinear operators on $C^\infty(M)$,
$C_0(f,g)=fg$, and $C_1(f,g)-C_1(g,f)=i\{f,g\}$.
A star product $\ast$ is called differential if all the 
operators $C_r, r\geq 0,$ are bidifferential.

Two differential star products $\ast$ and $\ast'$ are called 
equivalent if there exists an isomorphism of algebras
$B: (C^\infty(M)[[\nu]],\ast) \to (C^\infty(M)[[\nu]],\ast')$,
where $B = Id + \nu B_1 + \nu^2 B_2 +\dots$ and $B_r, r \geq 
1,$  are differential operators on $C^\infty(M)$.

If $(M,\omega)$ is a symplectic (and therefore a Poisson) 
manifold then to each differential star-product $\ast$ on $M$ 
one can relate its characteristic class $cl(\ast)\in 
(1/i\nu)[\omega]+H^2(M,\C)[[\nu]]$ (see 
\cite{Fedbook},\cite{NT1},\cite{NT2},\cite{Del}).
The equivalence classes of differential star-products on 
$(M,\omega)$ are bijectively parameterized by the elements of 
the affine vector space $(1/i\nu)[\omega]+H^2(M,\C)[[\nu]]$ 
via the mapping $\ast \mapsto cl(\ast)$.

The existence proofs and description of equivalence classes of 
star-products on symplectic manifolds were given by a number 
of  people (see \cite{DWL},\cite{OMY},\cite{Fed} for 
existence proofs and \cite{Fedbook},\cite{NT1},\cite{NT2},\cite{Del} 
for the classification). The questions of existence and 
classification on general Poisson manifolds were solved by 
Kontsevich in \cite{Kon}.

Historically the first examples of star products (Moyal-Weyl
and Wick star-products) were obtained from asymptotic
expansions of related symbol products w.r.t.  a numerical
parameter $\hbar$ (the `Planck constant'), as $\hbar\to 0$.
Later a number of star-products on K\"ahler manifolds were
obtained from the asymptotic expansion in $\hbar$ of the
symbol product of Berezin's covariant symbols (see
\cite{CGR2},\cite{CGR3},\cite{CGR4},\cite{M},\cite{MO},\cite{CMP2}).
These star-products on K\"ahler manifolds are
differential and have the property of `separation of
variables':  the corresponding bidifferential operators $C_r$
differentiate their first argument only in antiholomorphic
directions and the second argument only in holomorphic ones.
The deformation quantizations with separation of variables on
an arbitrary K\"ahler manifold where completely described and
parameterized by geometric objects, the formal deformations of 
the
K\"ahler (1,1)-form (see \cite{CMP1}).  The star-products
obtained from the product of Berezin's covariant symbols on
coadjoint orbits of compact semisimple Lie groups where
identified in \cite{CMP2} (that is, the corresponding formal
(1,1)-form was calculated).

Another interesting case where star-products are obtained from
symbol constructions is the so called Berezin-Toeplitz
quantization on arbitrary compact K\"ahler manifolds
\cite{BMS},\cite{Schldef},\cite{Schlbia95}.  In this
quantization scheme Berezin's contravariant symbols are used
(see \cite{Berequ}), which, in general, do not have a well
defined symbol product.  The semiclassical properties of
Berezin-Toeplitz quantization in \cite{BMS} were proved with
the use of generalized Toeplitz structures developed by
Guillemin and Boutet de Monvel in \cite{GBdM}.  With the same
technique for compact K\"ahler manifolds an associated
deformation quantization, the Berezin-Toeplitz deformation
quantization, was constructed \cite{Schldef},\cite{Schlbia95}.
However, this construction is very implicit and it is a
daunting task to calculate with its use the operators $C_r$
even for small values of $r$.  In \cite{KS} the
Berezin-Toeplitz deformation quantization was shown to be a
differentiable deformation quantization with separation of
variables and its classifying (1,1)-form was explicitly
calculated.  Thus the operators $C_r$ of the Berezin-Toeplitz
deformation quantization can be calculated recursively by the
simple algorithm from \cite{CMP1}.

It was shown by Guillemin in \cite{Guil} that, using
generalized Toeplitz structures, it is even possible to
construct a star-product on an arbitrary compact symplectic
manifold.  Later Borthwick and Uribe \cite{BU} introduced a
natural `almost-K\"ahler quantization' on an arbitrary compact
almost-K\"ahler manifold by defining an operator quantization
of Berezin-Toeplitz type.  From results by Guillemin and Uribe
\cite{GU} it follows that the generalized Toeplitz structure
associated to these operators exists and the correct
semi-classical limit in the almost K\"ahler case follows like
in \cite{BMS}.  Therefore there should exist the corresponding
natural deformation quantization on an arbitrary compact
almost-K\"ahler manifold. It would be interesting to describe 
this deformation quantization directly.

The goal of this Letter is to show that one can construct a
natural deformation quantization on an arbitrary
almost-K\"ahler manifold using Fedosov's machinery.  To this
end we use the natural affine connection introduced by
Yano \cite{Y}, which respects the almost-K\"ahler structure.
This connection necessarily has a nontrivial torsion whenever
the almost complex structure is non-integrable.  We generalize
Fedosov's construction to the case of affine connections with
torsion and obtain a star-product on an arbitrary
almost-K\"ahler manifold.  In the K\"ahler case this
star-product coincides with the star-product with separation
of variables (or of Wick type) constructed in \cite{BW}.  Our
considerations were very much motivated by the paper
\cite{DLS}, where the calculations concerning the star-product
from \cite{BW} were made in arbitrary local coordinates
(rather than in holomorphic coordinates as in \cite{BW}).

It follows from the results obtained in \cite{LMP1} and 
\cite{CMF} that the characteristic class of the star-product 
from \cite{BW} is $(1/i\nu)[\omega] -(1/2i)\varepsilon$, where 
$\varepsilon$ is the canonical class of the underlying 
K\"ahler manifold. The canonical class is, however, defined 
for any almost-complex manifold.

We calculate the crucial part of the characteristic class $cl$ 
of the star-product $\ast$ which we have constructed on an 
arbitrary almost-K\"ahler manifold, it's coefficient $c_0$ at 
the zeroth degree of $\nu$. As in the K\"ahler case, 
$c_0(\ast) = -(1/2i)\varepsilon$.

{\bf Acknowledgements.} We would like to thank 
B.~Fedosov and S.~Lya\-khovich for stimulating discussions, and the Referee for
useful comments.  A.K.
also thanks the DFG for financial support and the Department
of Mathematics and Computer Science of the University of
Mannheim, Germany, for their warm hospitality.

\section{A connection preserving an almost-K\"ahler 
structure}\label{S:con}
In this sections we recall some elementary properties of the 
Nijenhujs tensor and reproduce the construction from \cite{Y} 
of a natural affine connection on an arbitrary almost-K\"ahler 
manifold, which respects the almost-K\"ahler structure.
                                                                                                                   
Let $(M,J,\omega)$ be an almost-K\"ahler manifold, i.e., a
manifold $M$ endowed with an almost-complex structure $J$ and 
a symplectic form $\omega$
which are compatible in the following sense, 
$\omega(JX,JY)=\omega(X,Y)$ for any vector fields $X,Y$ on 
$M$, and $g(X,Y)=\omega(JX,Y)$ is a
Riemannian metric.  Actually, in what follows we only 
require the metric $g$ to be nondegenerate, but not 
necessarily positive definite. It is known that on any 
symplectic manifold $(M,\omega)$ one can choose a compatible 
almost complex structure to make it to an almost-K\"ahler 
manifold.

In local coordinates $\{ x^k\}$ on a chart $U\subset M$ set
\[
\p_j = \p/\p x^j,\ J(\p_j)=J_j^k \p_k,\
\omega_{jk}=\omega(\p_j,\p_k),\ g_{jk}=g(\p_j,\p_k)\ .
\]
Then
\begin{equation}\label{E:cs}
J_j^k=g_{j\alpha}\omega^{\alpha k} =
g^{k\alpha}\omega_{\alpha j},
\end{equation}
where $(g^{jk})$ and $(\omega^{jk})$ are inverse matrices for 
$(g_{kl})$ and $(\omega_{kl})$ respectively. Here, as well as 
in the sequel, we use the tensor rule of summation over 
repeated indices.

The Nijenhuis tensor for the complex structure $J$ is given by 
the following formula:
\begin{equation}\label{E:nij}
   N(X,Y) := [X,Y] + J[JX,Y] + J[X,JY] - [JX,JY],
\end{equation}
where $X,Y$ are vector fields on $M$. In local coordinates set 
$N(\p_j,\p_k) = N^l_{jk} \p_l$. Formula (\ref{E:nij}) takes 
the 
form
\begin{equation}\label{E:nijc}
   N^l_{jk} = (\p_\alpha J^l_j)J^\alpha_k - (\p_\alpha 
J^l_k)J^\alpha_j + (\p_j J^\beta_k - \p_k 
J^\beta_j)J^l_\beta.
\end{equation}
Let $X$ be a vector field of type (1,0) with respect to the
almost-complex structure $J$.  Then it follows from formula
(\ref{E:nij}) that 
\[
N(X,Y)=(1+iJ)([X,Y]+J[X,JY])\ .
\]  Therefore
$N(X,Y)$ is of type (0,1). The latter is also true if $Y$
is of type (1,0) instead.  Similarly, if $X$ or $Y$ is of
type (0,1) then $N(X,Y)$ is of type (1,0).  In
particular, if $X$ and $Y$ are of types (1,0) and (0,1)
respectively, then $N(X,Y)=0$.  Therefore, $(1\pm iJ)N((1\pm
iJ)\p_j,\p_k)=0$ and $N((1\pm iJ)\p_j,(1\mp iJ)\p_k)=0$,
from whence we obtain the two following formulas:
\begin{equation}\label{E:ident} 
N^l_{jk} = J^\alpha_j
N^\beta_{\alpha k} J^l_\beta \mbox{ and } N^l_{jk} = -
J^\alpha_j N^l_{\alpha\beta} J^\beta_k.  
\end{equation}
Let $\nabla$ be an arbitrary affine connection on $M$. In 
local 
coordinates let $\Gamma^l_{jk}$ be its Christoffel symbols and 
$T^l_{jk}=\Gamma^l_{jk}-\Gamma^l_{kj}$ the torsion tensor.

In the sequel we shall need the formula
\begin{equation}\label{E:cycl}
  \nabla_j \omega_{kl} + \nabla_k \omega_{lj} + \nabla_l 
\omega_{jk} = - ( T^\alpha _{jk}\omega_{\alpha l} + T^\alpha 
_{kl}\omega_{\alpha
j} + T^\alpha _{lj}\omega_{\alpha k}),
\end{equation}  
which easily follows from the closedness of $\omega$. Here 
$\nabla_i := \nabla_{\p_i}$.

By a direct calculation we get from formula (\ref{E:nijc}) 
that
\begin{equation}\label{E:nijtor}
   N^l_{jk} = (\nabla_\alpha J^l_j)J^\alpha_k - (\nabla_\alpha 
J^l_k)J^\alpha_j + (\nabla_j J^\beta_k - \nabla_k 
J^\beta_j)J^l_\beta - S^l_{jk},
\end{equation}
where
\begin{equation}\label{E:s}
S^l_{jk} = T^l_{jk} - J^\alpha_j T^l_{\alpha\beta} J^\beta_k + 
J^\alpha_j T^\beta_{\alpha k} J^l_\beta - J^\alpha_k 
T^\beta_{\alpha j} J^l_\beta.
\end{equation}
Assume that the connection $\nabla$ respects the 
almost-complex 
structure $J$. Then formula (\ref{E:nijtor}) reduces to 
\begin{equation}\label{E:red}
   N^l_{jk} = - S^l_{jk}.
\end{equation}
Thus (\ref{E:red}) is a necessary condition for $\nabla$ to 
respect the almost-complex structure $J$. 

One can check directly with the use of formulas 
(\ref{E:ident}) and 
(\ref{E:s})
that if
\begin{equation}\label{E:tors}
 T^l_{jk} = (-1/4)N^l_{jk}  
\end{equation}
then condition (\ref{E:red}) is always satisfied. 

\begin{proposition}\label{P:connect}
Let $\nabla$ be the unique affine connection 
which respects the metric $g$ and has the torsion given by 
formula (\ref{E:tors}). Then $\nabla$ also 
respects the symplectic form $\omega$ and therefore the 
complex 
structure $J$.
\end{proposition}

\noindent
This result is due to Yano \cite{Y}. For the 
convenience of the reader we provide here a proof of the 
statement.
\begin{proof}
Denote $Z_{jkl} = T^\alpha _{jk}\omega_{\alpha l} + T^\alpha 
_{kl}\omega_{\alpha
j} + T^\alpha _{lj}\omega_{\alpha k} = (-1/4)(N^\alpha 
_{jk}\omega_{\alpha l} + N^\alpha _{kl}\omega_{\alpha
j} + N^\alpha _{lj}\omega_{\alpha k}).$ Clearly, $Z_{jkl}$ is 
a 
totally antisymmetric tensor. 
It follows from  (\ref{E:cycl}) that
\begin{equation}\label{E:zcycl}
  \nabla_j \omega_{kl} + \nabla_k \omega_{lj} + \nabla_l 
\omega_{jk} = - Z_{jkl}.
\end{equation}  
First we show that
\begin{equation}\label{E:jz}
J^\alpha_j Z_{\alpha kl} = J^\alpha_k Z_{\alpha lj} = 
J^\alpha_l Z_{\alpha jk}.
\end{equation}
One has
\begin{equation}\label{E:intrm}
J^\alpha_j Z_{\alpha kl} = (-1/4)(J^\alpha_j N^\beta _{\alpha 
k} \omega_{\beta l} + J^\alpha_j N^\beta _{l\alpha} 
\omega_{\beta k} + J^\alpha_j N^\beta 
_{kl}\omega_{\beta\alpha}).
\end{equation}
Using (\ref{E:ident}) and (\ref{E:cs}) one obtains
$J^\alpha_j N^\beta _{\alpha k}\omega_{\beta l} = 
-N^\alpha_{jk}J^\beta_\alpha\omega_{\beta l} = -N^\alpha_{jk} 
g_{\alpha l}$. Similarly, $J^\alpha_j N^\beta _{l\alpha} 
\omega_{\beta k} = - N^\alpha _{lj} g_{\alpha k}.$
Finally, $J^\alpha_j N^\beta_{kl} \omega_{\beta\alpha} = - 
N^\alpha_{kl} g_{\alpha j}.$ Therefore
\[
J^\alpha_j Z_{\alpha kl} = (1/4)(N^\alpha_{jk} g_{\alpha l} + 
N^\alpha _{lj} g_{\alpha k} + N^\alpha_{kl} g_{\alpha j}),
\]
from whence (\ref{E:jz}) follows.

Formula (\ref{E:nijtor}) takes the form
\begin{equation}\label{E:eqn}
   (\nabla_\alpha J^l_j)J^\alpha_k - (\nabla_\alpha 
J^l_k)J^\alpha_j + (\nabla_j J^\beta_k - \nabla_k 
J^\beta_j)J^l_\beta = 0,
\end{equation}
due to the fact that (\ref{E:red}) is now satisfied.
Using (\ref{E:cs}) and the condition that $\nabla$ respects 
the metric $g$ we rewrite (\ref{E:eqn}) as follows:
\begin{equation}\label{E:eqnc}
 g^{l\beta}(\nabla_\alpha \omega_{\beta j})J^\alpha_k -
g^{l\beta}(\nabla_\alpha \omega_{\beta k})J^\alpha_j +
g^{\alpha\beta}(\nabla_j\omega_{\alpha k} - 
\nabla_k\omega_{\alpha j})J^l_\beta = 0.
\end{equation}  
Using (\ref{E:cs}) once more, we get 
$g^{\alpha\beta}J^l_\beta 
= \omega^{\alpha l} = - g^{l\beta} J^\alpha_\beta$ and
rewrite (\ref{E:eqnc}) in the form
\begin{equation}\label{E:eqnf}
 (\nabla_\alpha \omega_{\beta j})J^\alpha_k -
(\nabla_\alpha \omega_{\beta k})J^\alpha_j -
(\nabla_j\omega_{\alpha k} - \nabla_k\omega_{\alpha 
j})J^\alpha_\beta  = 0.
\end{equation}  
It follows from formula (\ref{E:zcycl}) that 
$\nabla_j\omega_{\alpha k} - \nabla_k\omega_{\alpha j} = 
\nabla_\alpha\omega_{jk} + Z_{\alpha jk}$. Therefore one gets 
from  (\ref{E:eqnf}) that
\begin{equation}\label{E:eqnjz}
 (\nabla_\alpha \omega_{\beta j})J^\alpha_k +
(\nabla_\alpha \omega_{k\beta})J^\alpha_j -
 (\nabla_\alpha \omega_{jk})J^\alpha_\beta = J^\alpha_\beta 
Z_{\alpha 
jk}. 
\end{equation}  
Cyclicly permuting the indices $\beta \to j \to k \to \beta$ 
in 
(\ref{E:eqnjz}) and adding the resulting equation to 
(\ref{E:eqnjz}) one obtains
\begin{equation}\label{E:fin}
 2(\nabla_\alpha \omega_{\beta j})J^\alpha_k =
 J^\alpha_\beta Z_{\alpha jk} + J^\alpha_j Z_{\alpha k\beta} = 
2J^\alpha_k Z_{\alpha \beta j}.
\end{equation}  
The last equality in (\ref{E:fin}) follows from (\ref{E:jz}). 
Thus
\begin{equation}\label{E:shrt}
\nabla_\alpha \omega_{\beta j} = Z_{\alpha \beta j}.
\end{equation}  
Summing up (\ref{E:shrt}) over all the cyclic permutations of 
indices $\alpha,\beta,j$, one gets from (\ref{E:zcycl}) that 
$-Z_{\alpha \beta j} = 3Z_{\alpha \beta j}$. Therefore 
$Z_{\alpha \beta j}=0$. 
Now the statement of the proposition follows from 
(\ref{E:shrt}) 
and (\ref{E:cs}).
\end{proof}

Notice that if $(M,J,\omega)$ is a K\"ahler manifold then the 
connection $\nabla$ from Proposition \ref{P:connect} is just 
the K\"ahler connection.

\section{A modification of Fedosov's construction}
\label{S:mod}
In this section we shall slightly modify Fedosov's 
construction 
to obtain a 
deformation quantization on an almost-K\"ahler manifold 
$(M,J,\omega)$ endowed with a 
fixed affine connection $\nabla$ which respects the 
almost-K\"ahler structure 
and has a nontrivial torsion. The existence of such a 
connection was shown in 
the previous section. 

As in Section \ref{S:con}, we shall work in local coordinates 
$\{x^k\}$ on a coordinate chart $U \subset M$ and use the same 
notation. Following Fedosov, denote by $\{y^k\}$ the fibre 
coordinates on the tangent bundle w.r.t. the frame $\{\p_k\}$.  

We introduce a tensor $\L^{jk} := \omega^{jk} - i g^{jk}$ on $M$ 
and define the  formal Wick algebra $W_x$ for $x\in M$ associated 
with the tangent space $T_x M$, whose elements are formal 
series
\[
    a(\nu,y) = \sum_{r \geq 0, |\alpha| \geq 0} \nu^r 
a_{r,\alpha} y^\alpha,
\]
where $\alpha$ is a multi-index and the standard multi-index 
notation is used. The formal Wick product on $W_x$ is given
by the formula    
\begin{equation}\label{E:wick}
    a \circ b\; (y) := \exp \left(\frac{i\nu}{2}\L^{jk} 
\frac{\p^2}{\p y^j\p z^k} \right) a(y)b(z)|_{z = y}.
\end{equation}
Taking a union of algebras $W_x$ we obtain a bundle $W$ of 
formal Wick algebras. Denote by $\W$ the sheaf of its smooth 
sections.
The fibre product (\ref{E:wick}) can be extended to the space 
$\W \otimes \L$ of $W$-valued differential forms by means of 
the usual exterior product of the scalar forms $\L$.

We introduce gradings $deg_\nu, deg_s, deg_a$ on 
$\W\otimes\L$ defined on homogeneous elements $\nu, y^k, dx^k$ 
as follows:
\[
deg_\nu(\nu)=1,\quad deg_s(y^k)=1,\quad deg_a(dx^k)=1\ .
\]
All other 
gradings of the elements $\nu, y^k, dx^k$ are set to zero. The 
grading $deg_a$ is induced from the standard grading on $\L$.

The product $\circ$ on $\W \otimes \L$ is bigraded w.r.t. the 
grading $Deg = 2 deg_\nu + deg_s$ and the grading $deg_a$.

The connection $\nabla$ can be extended to an operator on 
$\W \otimes \L$ such that for $a\in\W$ and a scalar 
differential form $\lambda$ 
\begin{equation}\label{E:conf}
\nabla (a \otimes \lambda) := \left( \frac{\p a}{\p x^j} - 
\Gamma^l_{jk} y^k \frac{\p a}{\p y^l}\right)\otimes 
(dx^j\wedge 
\lambda) + a \otimes d\lambda.
\end{equation}
Using formulas (\ref{E:wick}) and (\ref{E:conf}) one can show 
that $\nabla$ is a $deg_a$-graded derivation of the algebra
$(\W \otimes \L,\circ)$.

We introduce Fedosov's operators $\delta$ and $\delta^{-1}$ on 
$\W\otimes\L$ as follows. Assume $a\in \W\otimes\L$ is 
homogeneous w.r.t. the gradings $deg_s$ and $deg_a$ with 
$deg_s(a)=p,\ deg_a(a)=q$. Set
\begin{equation*}
   \delta(a) = dx^j \wedge \frac{\p a}{\p y^j} \mbox{\quad and 
\quad} \delta^{-1}a = 
\begin{cases}    
 \frac {1}{p+q}
 y^j i\left(\frac{\p}{\p x^j}\right) a 
   &\text{if  $p+q > 0$,}
\\
 0, &\text{if $p=q=0$.}
\end{cases}
\end{equation*}
Then for $a\in \W\otimes\L$ one has
\begin{equation}\label{E:hodge}
   a = \delta\delta^{-1}a + \delta^{-1}\delta a + \sigma(a),
\end{equation}
where $a \mapsto\sigma(a)$ is the projection on the 
$(deg_s,deg_a)$-bihomogeneous part 
of $a$ of bidegree zero ($deg_s(a)=deg_a(a)=0$).
It is easy to check that the operator $\delta$ is also a 
$deg_a$-graded derivation of the 
algebra $(\W \otimes \L,\circ)$.

Define the elements
\[
T:= \frac{1}{2}\, \omega_{s \alpha}T^\alpha_{kl} y^s dx^k 
\wedge 
dx^l
\mbox{\quad and \quad}
R := \frac{1}{4}\, \omega_{s\alpha} R^\alpha_{tkl}y^s y^t 
dx^k\wedge 
dx^l
\]
of $\W\otimes\L$, where 
\[
 R^s_{tkl} := \frac{\p \Gamma^s_{lt}}{\p x^k} - \frac{\p 
\Gamma^s_{kt}}{\p x^l} + \Gamma^s_{k \alpha}\Gamma^\alpha_{lt} 
- \Gamma^s_{l \alpha}\Gamma^\alpha_{kt}
\]
is the curvature 
tensor of the connection $\nabla$. Then the formulas 
\begin{equation}\label{E:ad}
    [\nabla,\delta]=\frac{i}{\nu}\, ad_{Wick}(T),\quad 
\nabla^2 
= -\frac{i}{\nu}\, ad_{Wick} (R)
\end{equation}
can be obtained by a direct calculation using (\ref{E:wick}),
(\ref{E:conf}) and the identity 
\begin{equation}\label{E:sari}
\omega_{s\alpha} R^\alpha_{tkl} = \omega_{t\alpha} 
R^\alpha_{skl}
\end{equation}
proved in \cite{GRS} (the proof is valid also for connections 
with torsion). In (\ref{E:ad}) $[\cdot,\cdot]$ is the 
$deg_a$-graded commutator of 
endomorphisms of $\W\otimes\L$ and $ad_{Wick}$ is defined via 
the $deg_a$-graded commutator in $(\W\otimes\L,\circ)$.

The following two theorems are minor modifications of the 
standard statements of Fedosov's theory adapted to the case of 
affine connections with torsion.
We shall denote the $Deg$-homogeneous component of degree $k$ 
of an element $a\in \W\otimes\L$ by $a^{(k)}$.
\begin{theorem}\label{T:fedcon}
There exists a unique element $r \in \W\otimes\L$ such that 
$r^{(0)}=r^{(1)}=0,\ deg_a(r)=1,\ \delta^{-1}r = 0$, 
satisfying the equation 
\[
\delta r = T + R + \nabla r - 
\frac{i}{\nu}\; r\circ r\ .
\]
It can be calculated recursively with respect to the total 
degree $Deg$ as follows:
\begin{gather*}
                 r^{(2)} = \delta^{-1} T,\\
   r^{(3)} = \delta^{-1}\left(R + \nabla r^{(2)} - 
   \frac{i}{\nu}\, r^{(2)}\circ r^{(2)} \right),\\
   r^{(k+3)} = \delta^{-1}\left(\nabla r^{(k+2)} - 
\frac{i}{\nu} \sum_{l=0}^{k} r^{(l+2)}\circ r^{(k-l+2)} 
\right), k\geq 
1.
\end{gather*}
 Then the Fedosov connection $D := -\delta + \nabla - 
\frac{i}{\nu} 
ad_{Wick} (r)$ is flat, i.e., $D^2=0$.
\end{theorem}
\noindent
The proof of the theorem is by induction, with the use of the 
identities
\[
\delta T = 0 \mbox{\quad and \quad} \delta R = \nabla T.
\]
The identity $\delta T = 0$ follows from (\ref{E:cycl}) 
and the fact that the connection $\nabla$ respects the form 
$\omega$. The identity $\delta R = \nabla T$ can be proved by 
a direct calculation with the use of (\ref{E:sari}).

The Fedosov connection $D$ is a $deg_a$-graded derivation of 
the algebra $(\W\otimes\L,\circ)$. Therefore $\W_D := \ker D 
\cap \W$ is a subalgebra of $(\W,\circ)$.
\begin{theorem}\label{T:fedquant}
The projection $\sigma: \W_D \to C^\infty(M)[[\nu]]$ onto the 
part of $deg_s$-degree zero is a bijection. The inverse 
mapping $\tau: C^\infty(M)[[\nu]] \to \W_D$ for a function $f 
\in C^\infty(M)$ can be calculated recursively w.r.t. the 
total degree $Deg$ as follows:
\begin{gather*}
      \tau(f)^{(0)} = f,\\
      \tau(f)^{(k+1)} = \delta^{-1} \left(\nabla \tau(f)^{(k)} 
- \frac{i}{\nu} \sum_{l=0}^{k} 
ad_{Wick}\bigl(r^{(l+2)}\bigr)\bigl(\tau(f)^{(k-l)}\bigr)
  \right), k\geq 0.
\end{gather*}
The product $\ast$ on $C^\infty(M)[[\nu]]$ defined by the 
formula 
\[
f \ast g := \sigma(\tau(f) \circ \tau(g))\ ,
\]
is a 
star-product on $(M,J,\omega)$.
\end{theorem}
Assume that the connection $\nabla$ in Theorem \ref{T:fedcon} 
is as in Proposition \ref{P:connect}.
Then in the case when $(M,J,\omega)$ is a K\"ahler manifold 
the star-product given by Theorem \ref{T:fedquant} coincides 
with the star-product of  Wick type constructed in 
\cite{BW}.

\section{Calculation of the class $c_0$}\label{S:class}
It is well known that to each star-product $\ast$ 
on a symplectic manifold $(M,\omega)$ a formal cohomology 
class $cl(\ast) \in (1/i\nu)[\omega] + H^2(M,\C)[[\nu]]$ is 
related (see, e.g., \cite{GR}). This class (named the 
characteristic class of deformation quantization) determines 
the star-product up to equivalence. Denote by $c_0(\ast)$ the 
coefficient of $cl(\ast)$ at zeroth degree of the formal 
parameter $\nu$. The class $c_0$ is, in some sense, the most 
intriguing part of the characteristic class $cl$. Only the 
coefficient $c_0(\ast)$ of the class 
$cl(\ast)=(1/i\nu)[\omega]+c_0(\ast)+\dots$ can not be 
recovered from Deligne's intrinsic class. 
Also the cohomology class of the formal K\"ahler form 
parameterizing a quantization with separation of variables on 
a K\"ahler manifold differs from the characteristic class of 
this quantization only in the coefficient $c_0$ (see 
\cite{LMP1}).

In this section we shall calculate the class $c_0$ of the 
deformation quantization obtained in Theorem \ref{T:fedquant}.
First we 
recall the definition of the class $c_0$ of a star-product 
(\ref{E:star}) (see, e.g., \cite{LMP1}). 
For a function $f\in C^\infty(M)$ on a symplectic manifold 
$(M,\omega)$ denote by $\xi_f$ the corresponding Hamiltonian 
vector field.
For a bilinear operator $C=C(f,g)$ denote by $C^-$ its 
antisymmetric part, $C^-(f,g):=(1/2)(C(f,g)-C(g,f))$.
A star-product $\ast$ given by (\ref{E:star}) is called {\it 
normalized} if 
$C_1(f,g)=(i/2)\{f,g\}$. For a normalized star-product  
$\ast$ the bilinear operator $C^-_2$ is a de Rham -- 
Chevalley 2-cocycle. There exists a unique closed 2-form 
$\varkappa$ such that for all $f,g\in C^\infty(M)$
 one obtains $C^-_2(f,g)= (1/2)\varkappa(\xi_f,\xi_g)$. The class $c_0$ of a 
normalized star-product $\ast$ is defined as  
$c_0(\ast) := [\varkappa]$.

It is well known that each star-product on a symplectic 
manifold is equivalent to a normalized one. One  defines the 
class $c_0(\ast)$ of a star-product $\ast$ as the cohomology class 
 $c_0(\ast')$ of an equivalent normalized star-product $\ast'$.
In order to calculate the class $c_0(\ast)$ of the star-product 
$\ast$ from Theorem \ref{T:fedquant} we shall first construct 
an equivalent normalized star-product $\ast'$.

We introduce the fibrewise equivalence operator on $\W$ 
defined by the formula
\begin{equation}\label{E:bertrans}
G := \exp \left(-\nu\Delta\right),
\end{equation}
where $\Delta$ is given in local coordinates as follows:
\[
     \Delta = \frac{1}{4} g^{jk}\frac{\p^2}{\p y^j y^k}.
\]
It is well known that the fibrewise star-product $\circ'$  
defined on $\W$ as follows, $a \circ' b := G(G^{-1} \circ 
G^{-1}b)$, is the Weyl star-product:
\begin{equation}\label{E:weyl}
       a \circ' b\; (y) = \exp \left(\frac{i\nu}{2}\omega^{jk} 
\frac{\p^2}{\p y^j\p z^k} \right) a(y)b(z)|_{z = y}.
\end{equation}
The following formulas 
\begin{equation}\label{E:commut}
   [\nabla,\Delta] = [\nabla, G] =0
\mbox{\quad and \quad} [\delta,\Delta] = [\delta, G] = 0
\end{equation}
can be checked directly.

Pushing forward the Fedosov connection $D$ obtained in Theorem 
\ref{T:fedcon} via $G$ and taking into account formulas 
(\ref{E:commut}) we obtain a connection 
\[
D' = GDG^{-1} = -\delta 
+ \nabla - \frac{i}{\nu}\, ad_{Weyl}(r')\ , 
\]
where $r' = Gr$ and 
$ad_{Weyl}$ is calculated with respect to the 
$\circ'$-commutator.

Denote by $\W_{D'} := \ker D' \cap \W$ the Fedosov subalgebra 
of the algebra $(\W,\circ')$.  Clearly, $\W_{D'} = GW_D$.
One can show just as in Theorem \ref{T:fedquant} that the 
restriction of the projection $\sigma$ to $\W_{D'},\ 
\sigma: \W_{D'} \to C^\infty(M)[[\nu]]$, is a bijection. 
Denote its inverse by $\tau'$. Then $f \ast' g := 
\sigma(\tau'(f)\circ' \tau'(g))$ is a star-product on 
$(M,\omega)$ which is equivalent to the star-product $\ast$. 
The operator $B: (C^\infty(M)[[\nu]],\ast) \to 
(C^\infty(M)[[\nu]],\ast')$ given by the formula $Bf = \sigma 
\bigl(G \tau (f)\bigr)$ establishes this equivalence.

{}From now on let $C_r, r\geq 1,$ denote the bidifferential 
operators defining the star-product $\ast'$. We have to show 
that $C_1(f,g) = (i/2)\{f,g\}$ and calculate $C^-_2$ in order 
to determine the class $c_0(\ast):=c_0(\ast')$.

For $a\in \W$ we prefer to write $a|_{y=0}$ instead of 
$\sigma(a)$.
For $f\in C^\infty(M)$ set 
$$\tau'(f) = t_0(f) + \nu t_1(f) 
+\nu^2 t_2(f)  +\dots\ .$$
 Since $\sigma(\tau'(f)) = 
\tau'(f)|_{y=0} = f$ we have $t_0(f)|_{y=0} = f$ and 
$t_r(f)|_{y=0} = 0$ for $r \geq 1$.
It follows from (\ref{E:weyl}) that for $f,g\in C^\infty(M)$
\begin{multline*}
f\ast' g = \left(\tau'(f)\circ' \tau'(g)\right)|_{y=0} =
\bigg(t_0(f)t_0(g) + 
\\
\nu \bigg(t_0(f) t_1(g) + t_1(f) t_0(g) + 
\frac{i}{2}\,\omega^{pq}\bigg( 
\frac{\p t_0(f)} {\p y^p} \frac{\p t_0(g)} {\p 
y^q}\bigg)\bigg) + \dots \bigg)|_{y=0} = 
\\
 fg + \frac{i\nu}{2}\,\omega^{pq}\bigg( \frac{\p t_0(f)} {\p 
y^p} \frac{\p t_0(g)} {\p y^q}\bigg)|_{y=0} + \dots,
\end{multline*}
from whence
\begin{equation}\label{E:cone}
   C_1(f,g) = \frac{i}{2}\,\omega^{pq}\left( \frac{\p t_0(f)} 
{\p y^p} \frac{\p t_0(g)} {\p y^q}\right)|_{y=0}.
\end{equation}

Similarly one obtains that
\begin{equation}\label{E:ctwo}
   C^-_2(f,g) = \frac{i}{2}\,\omega^{pq}\left( \frac{\p 
t_0(f)} {\p y^p} \frac{\p t_1(g)} {\p y^q} + \frac{\p t_1(f)} 
{\p y^p} \frac{\p t_0(g)} {\p y^q}\right)|_{y=0}.
\end{equation}
For an element $a^{(d)} \in \W\otimes\Lambda$ of $Deg$-degree $d$ 
denote by $a^{(d)}_s$ its homogeneous component of $deg_s$-degree $s$.
We have to calculate $\tau'(f)^{(1)}$ and the component 
$\tau'(f)^{(3)}_1$ of $\tau'(f)^{(3)}$. 
It follows from the condition $D' \tau'(f)=0$ that
\begin{equation}\label{E:lev}
 \tau'(f)^{(k+1)} = \delta^{-1} \left(\nabla \tau'(f)^{(k)} 
- \frac{i}{\nu} \sum_{l=0}^{k} 
ad_{Weyl}\bigl((r')^{(l+2)}\bigr)\bigl(\tau'(f)^{(k-l)}\bigr)
  \right), k\geq 0.
 \end{equation}
Since $\tau'(f)^{(0)} = f$, we get from (\ref{E:lev}) for 
$k=0$ that
$\tau'(f)^{(1)} = (\p f/\p x^p) y^p$. Therefore $(\p t_0(f)/\p 
y^p)|_{y=0} = \p f/\p x^p$.
Now (\ref{E:cone}) 
implies that $C_1(f,g) = (i/2)\{f,g\}$, i.e., that the 
star-product $\ast'$ is normalized.

Taking into account that $\tau'(f)^{(2)}$ is of $deg_s$-degree 
2, we obtain from (\ref{E:lev}) for $k=2$ that 
\begin{equation}\label{E:taup}
\tau'(f)^{(3)}_1 = - \frac{i}{\nu} \delta^{-1}\left(
ad_{Weyl}\bigl((r')^{(3)}_1\bigr)\bigl(\tau'(f)^{(1)}\bigr)
\right) = \delta^{-1}\left(\omega^{pq}\frac{\p 
\bigl((r')^{(3)}_1\bigr)}{\p y^p}\frac{\p f}{\p x^q}\right).
\end{equation}
Before calculating $(r')^{(3)}_1$ one can directly derive from 
(\ref{E:ctwo}) and (\ref{E:taup}) the following formula:
\begin{equation}\label{E:kapp}
   \varkappa = \frac{i}{\nu}\, 
   \delta\bigl((r')^{(3)}_1 \bigr).
\end{equation}
Denote by $[\cdot,\cdot]_\circ$ the commutator with respect to the 
Wick multiplication $\circ$. 
We shall need the following technical lemma
which can be proved by a straightforward calculation.
\begin{lemma}\label{L:aux}
Let $a = a^{(2)}_2, b = b^{(2)}_2$ be two homogeneous elements of $\W$, 
$(1/\nu)[a,b]_\circ = c^{(2)} = c^{(2)}_0 + c^{(2)}_2$, then 
$c^{(2)}_0 = \nu \Delta \left(c^{(2)}_2 \right)$.
\end{lemma}
Using the fact that the operator $G$ respects the total 
grading $Deg$ and $\delta$ lowers both $Deg$- and 
$deg_s$-gradings by 1, one can obtain from formula 
(\ref{E:commut}), and the formula $r' = Gr$ that
\begin{equation}\label{E:frst}
   \delta\bigl((r')^{(3)}_1 \bigr) = \delta r^{(3)}_1 - 
   \nu\Delta \left(\delta r^{(3)}_3 \right).
\end{equation}
We get from Theorem \ref{T:fedcon} that 
\begin{equation}\label{E:sec}
      \delta r^{(3)} = R + \nabla r^{(2)} - \frac{i}{\nu}\, r^{(2)}\circ r^{(2)},
\end{equation} 
where
\[
           r^{(2)} = \delta^{-1} T = \frac{1}{3}\, \omega_{s\alpha} 
T^\alpha_{tl} y^s y^t dx^l.
\]
Since the element $r^{(2)}$ is of $deg_a$-degree 1, we have 
\begin{equation}\label{E:thrd}
      \frac{i}{\nu}\, r^{(2)}\circ r^{(2)} = \frac{i}{2 \nu}\, [r^{(2)}, r^{(2)}]_\circ 
          = c^{(2)} = c^{(2)}_0 + c^{(2)}_2.
\end{equation}
We obtain from  (\ref{E:sec}) and (\ref{E:thrd}) that
\begin{equation}\label{E:four}
        \delta r^{(3)}_1   = - c^{(2)}_0  \mbox{\quad and \quad}
      \delta r^{(3)}_3 = R + \nabla r^{(2)} - c^{(2)}_2.
\end{equation} 
It follows from (\ref{E:kapp}), (\ref{E:frst}),  (\ref{E:four}), and Lemma \ref{L:aux}
that  
\begin{equation}\label{E:kappa}
   \varkappa = - i \Delta
\bigl(R + \nabla r^{(2)}\bigr) = - \frac{i}{8} J^t_s 
R^s_{tkl} dx^k \wedge dx^l - i \lambda,
\end{equation}   
where $\lambda = \Delta\bigl(\nabla r^{(2)}\bigr)$.
Introduce a global differential one-form $\mu = (1/6) J^t_s 
T^s _{tl} dx^l$ on $M$. A direct calculation shows that 
$\lambda = d \mu$, therefore the form $\lambda$ is exact.

Recall the definition of the canonical class $\varepsilon$ of 
an almost complex manifold $(M,J)$. The class $\varepsilon$ is 
the first Chern class of the subbundle $T'_\C M$ of vectors of 
type (1,0) of the complexified tangent bundle $T_\C M$. To 
calculate the canonical class of the almost-K\"ahler manifold 
$(M,J,\omega)$ take the same affine connection $\nabla$ on $M$ 
as that used in the construction of the star-product $\ast$.
Denote by $\hat R = (1/2) R^s_{tkl} dx^k \wedge dx^l$ 
the curvature matrix of the connection 
$\nabla$ and by $\Pi = (1/2)(Id - iJ)$ the projection operator 
onto the (1,0)-subspace. It follows immediately from 
(\ref{E:sari}) that $R^t_{tkl} = 0$ (see \cite{GRS}), i.e.,
$Tr \hat R =0$. The matrix $\Pi \hat R \Pi$ is the curvature 
matrix of the restriction of the connection $\nabla$ to 
$T'_\C M$. The Chern-Weyl form 
\[
\gamma = (1/i)Tr(\Pi \hat R 
\Pi) = (1/i)Tr (\Pi \hat R) = (-1/4)J^t_s R^s_{tkl} dx^k 
\wedge dx^l
\]
is closed. The canonical class is, by definition, 
$\varepsilon :=[\gamma]$. Now it is clear from (\ref{E:kappa}) 
that 
\begin{equation}
c_0(\ast) = [\varkappa] = -(1/2i)\varepsilon.
\end{equation}

\end{document}